\documentclass[graybox, envcountchap]{./svmult}

\usepackage{graphicx}           
\usepackage{makeidx}            
\usepackage{multirow}
\usepackage{multicol}           
\usepackage{caption}
\usepackage[bottom]{footmisc}   

\usepackage[utf8]{inputenc}
\usepackage[T1]{fontenc}
\usepackage{float}
\usepackage{wrapfig}
\usepackage{amsmath}
\usepackage{amsopn}
\usepackage{amsfonts}
\usepackage{hyperref}
\usepackage{lmodern}
\usepackage[english]{babel}
\usepackage{booktabs}  
\usepackage{enumitem}  
\usepackage{tabularx}
\usepackage[normalem]{ulem}
\usepackage{enumitem}

\usepackage{listings}
\usepackage{xcolor}

\usepackage{tikz}
\usetikzlibrary{arrows.meta,decorations.markings}
\usetikzlibrary{decorations.pathmorphing,decorations.pathreplacing,decorations.shapes}

\definecolor{royalpurple}{rgb}{0.47, 0.32, 0.66}

\definecolor{SegA}{RGB}{0,114,178}     
\definecolor{SegB}{RGB}{230,159,0}     
\definecolor{SegC}{RGB}{0,158,115}     

\usepackage{hyperref}
\usepackage{doi}

\newcommand{\R}{\mathbb{R}}

\usepackage{epstopdf}
\ifpdf
  \DeclareGraphicsExtensions{.eps,.pdf,.png,.jpg}
\else
  \DeclareGraphicsExtensions{.eps}
\fi

\usepackage{ifthen}
\newcommand{\mycomment}[1]{%
  \ifthenelse{\isodd{\value{page}}}{%
    \normalmarginpar%
    \marginpar{\tiny {#1}}%
  }{%
    \reversemarginpar%
    \marginpar{\tiny {#1}}%
  }}%

\newcommand{\TheTitle}{Overlapping Schwarz Preconditioners for  Pose-Graph SLAM in Robotics}

\newcommand{\TheShortTitle}{Overlapping Schwarz in Robotics}

\newcommand{\TheAuthors}{Stephan~Köhler$^1$, Oliver~Rheinbach$^{1,2}$, Yue~Xiang~Tee$^1$, and Sebastian~Zug$^1$}

\ifpdf
\hypersetup{
  pdftitle={\TheTitle},
  pdfauthor={\TheAuthors}
}
\fi

\usepackage{amsmath}
\usepackage{xcolor}
\usetikzlibrary{calc}

\definecolor{SegA}{RGB}{50,120,200}
\definecolor{SegB}{RGB}{220,120,40}
\definecolor{Common}{RGB}{90,90,90}

\tikzset{
  baseedge/.style={black, line width=0.9pt, line cap=round},
  subedgeA/.style={SegA, line width=2.0pt, line cap=round},
  subedgeB/.style={SegB, line width=2.0pt, line cap=round},
  nodebase/.style={circle, fill=black, draw=white, line width=0.4pt, inner sep=0pt, minimum size=6.5pt},
  nodeA/.style={circle, fill=SegA, draw=white, line width=0.4pt, inner sep=0pt, minimum size=6.5pt},
  nodeB/.style={circle, fill=SegB, draw=white, line width=0.4pt, inner sep=0pt, minimum size=6.5pt},
  ringA/.style={circle, fill=none, draw=SegA, line width=1.1pt, inner sep=0pt, minimum size=13pt},
  ringB/.style={circle, fill=none, draw=SegB, line width=1.1pt, inner sep=0pt, minimum size=13pt},
  conn/.style={dashed, black!35, line width=0.5pt},
  label/.style={font=\footnotesize},
  ilabel/.style={font=\small}
}

\newcommand{\basechain}{%
  \foreach \i in {1,...,9} {\coordinate (p\i) at ({\i-1},0);}
  \draw[baseedge] (p1) -- (p2) -- (p3) -- (p4) -- (p5) -- (p6) -- (p7) -- (p8) -- (p9);
  \foreach \i in {1,...,9} {
    \node[nodebase] at (p\i) {};
    \node[label, below=4pt] at (p\i) {$p_{\i}$};
  }
}

\begin{document}

\titlerunning{\TheShortTitle}
\title*{{\TheTitle}}
\author{\TheAuthors}
\institute{
$^1$ Faculty of Mathematics and Computer Science,
Technische Universität Bergakademie Freiberg,
Akademiestr.~6, 09599~Freiberg.\\
$^2$ Universitätsrechenzentrum (URZ),
Technische Universität Bergakademie Freiberg,
Bernhard-von-Cotta-Str. 1.\\
\email{stephan.koehler@math.tu-freiberg.de, oliver.rheinbach@math.tu-freiberg.de, Yue-Xiang.Tee@math.tu-freiberg.de, Sebastian.Zug@informatik.tu-freiberg.de}, 
}

\numberwithin{theorem}{section}
\numberwithin{equation}{section}

\maketitle

\abstract{
We investigate scalable two-level overlapping Schwarz domain decomposition methods with energy-minimizing coarse spaces of GDSW type (Generalized Dryja--Smith--Widlund type) as preconditioners for the 
sparse linear systems arising in graph-based nonlinear least-squares problems, specifically the pose-graph optimization back-end in Simultaneous Localization and Mapping (SLAM).
After a brief introduction to SLAM and domain decomposition preconditioners, we describe the nonlinear least-squares formulation, its linearization, and the resulting matrix structure,
to facilitate access
for readers without prior knowledge of either field. Numerical experiments demonstrate the numerical scalability of the preconditioned conjugate gradient method (CG): Using the two-level overlapping Schwarz preconditioner, the number of CG iterations remains bounded independently of the problem size, overcoming the typical limitations of simple preconditioners, including one-level Schwarz approaches. We further show that a simplified SLAM problem can be interpreted as a finite element problem using linear elastic bars, reinforcing the analogy to continuum mechanics and motivating the use of scalable domain decomposition techniques.
}



\section{Introduction}
Simultaneous Localization and Mapping (SLAM) addresses the problem of estimating a robot’s trajectory while simultaneously constructing a map of its environment. In large-scale scenarios such as long-term operation in outdoor environments, SLAM leads to very large sparse 
minimization problems derived from noisy sensor measurements. Solving these problems efficiently and robustly remains a major challenge, motivating the development of fast SLAM algorithms.



Preconditioning techniques have been proposed to accelerate the iterative solution of the large sparse linear systems arising in graph-based SLAM. 
These systems typically arise as the normal equations of the linearized subproblems solved at each iteration of nonlinear least-squares methods such as Gauss--Newton or Levenberg--Marquardt.

Different approaches have been explored in the literature.
For instance, a subgraph-based preconditioner constructed from a sparse subgraph structure, such as a spanning tree, was proposed in~\cite{dellaert2010subgraph} and later extended to incremental settings suitable for online SLAM in~\cite{ispcg_2014}.
Block-Jacobi preconditioning and the direct solution of Schur complement
systems are considered in~\cite{blockjacobi_2010}.
Incomplete Cholesky factorizations have also been used as preconditioners
for conjugate gradients applied to the sparse linear systems arising in
Levenberg--Marquardt SLAM~\cite{spa_2010}, where they are also compared to sparse direct solvers.
In practical SLAM systems such as the g2o software library~\cite{g2o}, the linear systems arising in pose graph optimization are typically solved using sparse direct Cholesky factorization (e.g., CHOLMOD).
While iterative methods such as preconditioned conjugate gradients are also available, they typically employ computationally inexpensive preconditioners such as block-Jacobi or diagonal scaling.

However, numerical scalability with growing problem size remains a challenge for many of these approaches, i.e., their performance tends to deteriorate as the problem size increases.
As emphasized in the overview paper~\cite{pastpresentfutureslam_2016}, the development of scalable SLAM methods remains a fundamental open challenge for achieving long-term autonomy in large environments.

In this paper, we demonstrate numerically that domain decomposition preconditioners, specifically two-level additive overlapping Schwarz preconditioners, applied to a simple test problem in SLAM, can lead to a numerically scalable method, i.e., the number of preconditioned conjugate gradient iterations remains small and bounded independently of the problem size.
Keeping the number of iterations small is important in practice, since large iteration counts may lead to loss of orthogonality in the conjugate gradient method in finite precision arithmetic.



The present version substantially revises the previous version of this manuscript. In particular, the numerical results for the one-level Schwarz method have been corrected; they were affected by an implementation error that introduced an unintended global coupling. As a consequence, the previous version incorrectly suggested scalability of the one-level method. The corrected experiments show that scalability is obtained only with the two-level Schwarz method; see Table~\ref{tab:weak_numerical}. In addition, the one-dimensional illustrative problem has been moved to the appendix in order to focus the main text on the two-dimensional nonlinear SLAM model problem.

We first give a brief introduction to both SLAM and overlapping Schwarz methods. This is followed by a description of the mathematical formulation of the problem and the problem setting. Finally, we present and discuss convergence results for one-level and two-level overlapping Schwarz preconditioners applied to a two-dimensional SLAM model problem. In the appendix, we also present a simple one-dimensional SLAM problem to illustrate connections between SLAM and continuum mechanical problems discretized by the finite element method.


\section{Simultaneous Localization and Mapping (SLAM)}
In SLAM, the trajectory of a robot and a map of the environment are estimated simultaneously from noisy sensor measurements.
%
%
The localization problem has been investigated since the 1990s, initially with a strong focus on filter-based approaches such as the extended Kalman filter; see, e.g.,~\cite{durrantwhyte2006slam1, durrantwhyte2006slam2, 10.1007/978-1-4471-1021-7_69,alsinet2008slam}. With the growing interest in SLAM in recent years, several comprehensive surveys have been published. In particular,~\cite{pastpresentfutureslam_2016} provides an 
overview of the development of SLAM, while more recent surveys discuss the architecture of complete SLAM systems and modern sensor configurations~\cite{alsadik2021simultaneous,chen2022slam}.


\begin{figure}[th]
  \centering
  \includegraphics[width=\linewidth]{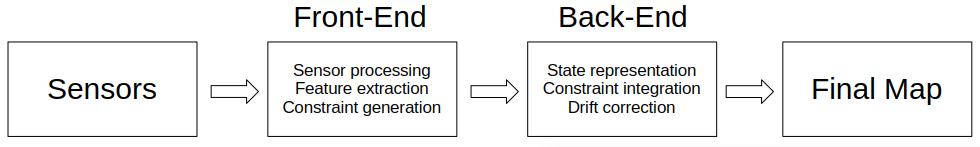}
  \caption{Simplified flow chart of a complete SLAM system.}\label{fig:slam_overview}
\end{figure}

A SLAM system starts with the collection of data from sensors, and ends with the creation of a map, with which the robot is able to orient itself. This process can be briefly summarized as in Fig.~\ref{fig:slam_overview}.

A SLAM system typically consists of a front-end and a back-end. The front-end processes sensor data to estimate relative motion (odometry) between robot poses. The specific approach depends on the sensor modality, for example visual SLAM for camera-based systems~\cite{macario2022comprehensive} or LiDAR-based SLAM for laser sensors~\cite{huang2021review}. By extracting features and performing data association, the front-end produces measurement constraints that are passed to the back-end for global estimation.

Due to the inevitable accumulation of noise in sensor measurements, drift increases over time. When the system revisits a previously observed location, additional constraints known as \emph{loop closures} can be introduced, allowing the accumulated drift to be corrected.


The back-end estimates the robot states from the motion and measurement constraints.
This can be achieved using either filter-based or graph-based methods.
Filter-based methods estimate the current state of the system recursively by propagating and updating a probability distribution as new measurements arrive.
In graph-based methods, robot poses, map variables, and measurements are represented as nodes connected by constraints in a graph, and the most consistent global solution is computed by solving a nonlinear least-squares problem; see~\cite{grisetti2011tutorial} for a tutorial.
Efficient implementations of such formulations are provided by general graph-optimization frameworks such as g2o~\cite{g2o}.


In this paper, after linearization 
using Gauss--Newton, the ill-condit\-ioned linear system is solved by a preconditioned conjugate gradient method using an additive overlapping Schwarz preconditioner.
This is described in more detail in the following sections.

\subsection{Mathematical Formulation of Pose Graph SLAM in Two Dimensions}\label{sec:mathslam}
We describe the specific nonlinear least-squares problem underlying the numerical experiments formulated in classical terms of finite-dimensional nonlinear optimization.

{\color{black}In two-dimensional SLAM, a {pose} represents the configuration of the robot in the plane. It consists of a position vector and an orientation angle and is written as
\[
p = (x,y,\theta)^T \in \mathbb{R}^2 \times (-\pi,\pi],
\]
where the orientation angle \(\theta\) is defined modulo \(2\pi\).
%
Given two poses
\[
p_i = (x_i,y_i,\theta_i)
\quad\text{and}\quad
p_j = (x_j,y_j,\theta_j),
\]
their {relative pose} is the configuration of \(p_j\) expressed in the local coordinate system of the previous pose \(p_i\).

We model the SLAM problem by a pose graph
$
G=(V,E),
$
where the vertex set
\[
V=\{0,1,\dots,N\}
\]
indexes the robot poses \(p_i\), and the set of directed edges
\[
E \subset V \times V
\]
contains relative pose measurements between pairs of poses. The edges include 
odometry measurements between consecutive poses as well as loop closures between poses which are revisited:
} %
%
An edge \((i,j) \in E\) exists in the pose graph if \(p_i\) and \(p_j\) are consecutive poses, i.e., \(j = i+1\), or if the edge $(i,j)$ represents a loop closure constraint between \(p_i\) and \(p_j\); see Fig.~\ref{fig:slam-loop-closure}.
\begin{figure}[tbh]
  \centering
  \begin{tikzpicture}[
    >=Stealth,
    truthnode/.style={
      circle, fill=red!80, draw=red!60!black,
      inner sep=0pt, minimum size=7pt
    },
    measnode/.style={
      circle, fill=blue!70, draw=blue!40!black,
      inner sep=0pt, minimum size=7pt
    },
    truthedge/.style={red!90!black, line width=1.2pt},
    measedge/.style={blue!60!black, line width=1.0pt, dashed},
    ]

    \def\mAx{ 0.045}  \def\mAy{-0.038}
    \def\mBx{ 1.038}  \def\mBy{-0.033}
    \def\mCx{ 1.041}  \def\mCy{ 1.028}
    \def\mDx{-0.032}  \def\mDy{ 1.052}
    \def\mEx{-0.112}  \def\mEy{ 0.018}
    \def\mFx{ 0.922}  \def\mFy{ 0.063}

    \begin{scope}[scale=3.2]

      \draw[truthedge] (0,0) -- (1,0) -- (1,1) -- (0,1) -- (0,0);
      \node[truthnode] (gt0) at (0,0) {};
      \node[truthnode] (gt1) at (1,0) {};
      \node[truthnode] (gt2) at (1,1) {};
      \node[truthnode] (gt3) at (0,1) {};

      \draw[measedge] (\mAx,\mAy) -- (\mBx,\mBy)
      -- (\mCx,\mCy)
      -- (\mDx,\mDy)
      -- (\mEx,\mEy)
      -- (\mFx,\mFy);

      \draw[green!80!black, line width=0.6pt, decorate, decoration={
        zigzag,
        segment length=2pt,
        amplitude=1pt,
        pre=lineto,  pre length=4.5pt,
        post=lineto, post length=4.5pt
      }]
      (\mEx,\mEy) -- (\mAx,\mAy);

      \node[measnode] (m0) at (\mAx,\mAy) {};
      \node[measnode] (m1) at (\mBx,\mBy) {};
      \node[measnode] (m2) at (\mCx,\mCy) {};
      \node[measnode] (m3) at (\mDx,\mDy) {};
      \node[measnode] (m4) at (\mEx,\mEy) {};
      \node[measnode] (m5) at (\mFx,\mFy) {};

      \node[red!80!black, font=\tiny, anchor=south west] at (gt0) {\,$p_0^*$};
      \node[red!80!black, font=\tiny, anchor=south west] at (gt1) {\,$p_1^*$};
      \node[red!80!black, font=\tiny, anchor=north east] at (gt2) {$p_2^*$\,};
      \node[red!80!black, font=\tiny, anchor=north west] at (gt3) {\,$p_3^*$};

      \node[blue!70!black, font=\tiny, anchor=north east] at (m0) {\,$\tilde{p}_0$};
      \node[blue!70!black, font=\tiny, anchor=north west] at (m1) {\,$\tilde{p}_1$};
      \node[blue!70!black, font=\tiny, anchor=south west] at (m2) {$\tilde{p}_2$};
      \node[blue!70!black, font=\tiny, anchor=south east] at (m3) {\,$\tilde{p}_3$};
      \node[blue!70!black, font=\tiny, anchor=south east] at (m4) {$\tilde{p}_4$\,};
      \node[blue!70!black, font=\tiny, anchor=south east] at (m5) {$\tilde{p}_5$\,};

      \node[truthnode,
      label={[red!70!black, font=\scriptsize]right:~ground truth $p^*$}]
      at (0.05, -0.22) {};
      \node[measnode,
      label={[blue!70!black, font=\scriptsize]right:~odometry $\tilde{p}$}]
      at (0.05, -0.33) {};
      \draw[green!80!black, font=\scriptsize, decorate, decoration={
        zigzag,
        segment length=2pt,
        amplitude=1pt,
        pre=lineto,  pre length=1pt,
        post=lineto, post length=1pt
      }] (0.0, -0.44) -- (0.1, -0.44);
      \node[green!80!black, font=\scriptsize] at (0.38, -0.44) {loop closure};

    \end{scope}  
%
    \def\bxl{-0.84}  \def\bxr{0.51}
    \def\byb{-0.48}  \def\byt{0.48}

    \draw[black!70, line width=0.9pt, dashed, rounded corners=2pt] (\bxl,\byb) rectangle (\bxr,\byt);

    \def\ixl{5.00}  \def\ixr{7.60}
    \def\iyb{1.00}  \def\iyt{3.00}

    \draw[black!40, line width=0.6pt] (\bxr,\byt) -- (\ixl,\iyt);
    \draw[black!40, line width=0.6pt] (\bxr,\byb) -- (\ixl,\iyb);

    \begin{scope}
      \clip (\ixl,\iyb) rectangle (\ixr,\iyt);

      \begin{scope}[xshift=6.60cm, yshift=1.96cm, scale=8,
        truthnode/.style={
          circle, fill=red!80, draw=red!60!black,
          inner sep=0pt, minimum size=6pt},
        measnode/.style={
          circle, fill=blue!70, draw=blue!40!black,
          inner sep=0pt, minimum size=6pt},
        ]
        \draw[green!80!black, line width=0.5pt, decorate, decoration={
          zigzag,
          segment length=2pt,
          amplitude=1pt,
          pre=lineto,  pre length=5pt,
          post=lineto, post length=5pt
        }](\mEx, \mEy) -- (\mAx,\mAy);

        \draw[measedge] (\mAx,\mAy) -- (\mBx,\mBy); 
        \draw[measedge] (\mDx,\mDy)-- (\mEx,\mEy) -- (\mFx,\mFy); 

        \draw[truthedge] (0,0) -- (1,0);
        \draw[truthedge] (0,1) -- (0,0);
        \node[truthnode] at (0,0) {};

        \node[measnode] (zi0) at ( 0.045,-0.038) {};   
        \node[measnode] (zi4) at (-0.112, 0.018) {};   

        \node[red!80!black,  font=\scriptsize, anchor=south west] at (0,0)          {$p_0^*$};
        \node[blue!70!black, font=\scriptsize, anchor=north east] at (\mAx,\mAy) {$\tilde{p}_0$};
        \node[blue!70!black, font=\scriptsize, anchor=south east] at (\mDx,\mDy) {$\tilde{p}_4$};
      \end{scope}
    \end{scope}

    \draw[black!70, line width=0.9pt, rounded corners=2pt] (\ixl,\iyb) rectangle (\ixr,\iyt);

    \node[font=\scriptsize\itshape, black!60, anchor=north] at ({(\ixl+\ixr)/2}, \iyb) {zoom of loop closure};
  \end{tikzpicture}
  \caption{Example for a SLAM problem; ground truth (red) and odometry data (blue); between the pose $\tilde{p}_{4}$ and $\tilde{p}_{0}$ is a loop closure edge (green)}\label{fig:slam-loop-closure}
\end{figure}


Algebraically, the relative pose is the pose of \(p_j\) expressed in the
local coordinate system of \(p_i\). It is given by the smooth nonlinear
mapping \(\Phi : \mathbb{R}^3 \times \mathbb{R}^3 \to \mathbb{R}^3\),
\begin{equation}
  \label{eq:slam-relative-pose}
  \Phi(p_i, p_j) = p_i^{-1} \circ p_j =
  \begin{pmatrix}
    R(\theta_i)^T (t_j - t_i) \\
    \theta_j - \theta_i
  \end{pmatrix},
\end{equation}
where \(t_i = (x_i,y_i)^T\) and
\begin{equation*}
  R(\theta_i) =
  \begin{pmatrix}
    \cos \theta_{i} & -\sin \theta_{i} \\
    \sin \theta_{i} & \cos \theta_{i}
  \end{pmatrix}
\end{equation*}
denotes the rotation matrix associated with \(\theta_i\). Multiplication
by \(R(\theta_i)^T\) transforms the global displacement \(t_j-t_i\) into the
local coordinate frame of pose \(p_i\). The angle difference is again mapped
to a fixed interval.

The problem data consist of relative pose measurements.  For certain index pairs $(i,j)$, a measurement
\begin{equation}
  \widetilde{\Phi}_{ij} = (\delta\tilde{x}_{ij},\, \delta\tilde{y}_{ij},\, \delta\tilde{\theta}_{ij})^T
\end{equation}
%
is given, represented by the edge $(i,j)\in E$ of the pose graph.
This edge represents a (potentially noisy) observation of the relative pose between $p_i$ and $p_j$, i.e., $\widetilde{\Phi}_{ij}\approx\Phi(p_i,p_j)$.
In our numerical experiments, two types of measurements are used:
\begin{itemize}
\item \emph{Sequential measurements}, which connect consecutive poses and incremental motion.
\item \emph{Loop closure measurements}, which connect poses that are far apart in the sequence and introduce global coupling.
\end{itemize}
Loop closures occur whenever the robot recognizes a place it has already visited. This allows for the correction of drift accumulated on the robot's path.
A loop closure is simply an additional constraint between poses indicating that these poses should be the same; see Fig.~\ref{fig:slam-loop-closure}. These constraints are treated in the same way as the other measurements, i.e., they participate in the least-squares adjustment and are not enforced exactly.
However, each measurement is weighted using a
diagonal (or symmetric positive definite) weight matrix $W$ with entries
\begin{equation}
W_{ij} \in \mathbb{R}^{3 \times 3},
\label{eq:weight}
\end{equation}
and a higher weight is often used for the loop closure measurements, indicating their higher importance.
Let us remark that we can set the entry for the orientation in the weight matrix to zero if we are not interested in the robot's orientation at a loop closure.

For a given estimate $\widetilde{\Phi}_{ij}$ for the relative pose $\Phi(p_i,p_j)$, the corresponding residual is defined as
\begin{equation}
  r_{ij}(p_i,p_j) = \Phi(p_i,p_j) - \widetilde{\Phi}_{ij}.
\label{eq:slam-residual-small}
\end{equation}
The rotational component of the residual is wrapped to the principal interval $(-\pi,\pi]$, thereby avoiding artificial angle errors differing by multiples of $2\pi$.
The residual measures the discrepancy between the predicted relative pose $\Phi(p_i,p_j)$ and the measured relative pose $\widetilde{\Phi}_{ij}$.


\subsection{Least-Squares Formulation}
\label{seq:lsq}
Using the residual defined in~(\ref{eq:slam-residual-small})
the
global least squares objective function is given by
\begin{equation}
  \mathcal{J}(p_1,\dots,p_N) = \frac{1}{2}\sum_{(i,j)\in E} {r_{ij}(p_i,p_j)}^T\, W_{ij}\, r_{ij}(p_i,p_j),
  \label{eq:slam-obj}
\end{equation}
where $W_{ij} \in \mathbb{R}^{3 \times 3}$
is the weight matrix introduced in (\ref{eq:weight}). 

The objective function is invariant under global rigid transformations.  The nonlinear least-squares problem is given by
\begin{equation}
  \label{eq:slam-min}
  \min_{(p_1,\dots,p_N)\in\R^{3N}}\mathcal{J}(p_1,\dots,p_N).
\end{equation}
In our benchmark problem, the null space of the problem is removed by fixing the initial pose.

Note that~\eqref{eq:slam-min} can also be written as
\begin{equation*}
  \min_{(p_1,\dots,p_N)\in\R^{3N}}\frac{1}{2}\sum_{(i,j)\in E_{\text{odom}}} {\|r_{ij}(p_i,p_j)\|}^2_{W_{ij}} + \frac{1}{2}\sum_{(i,j)\in E_{\text{loop}}}{\|r_{ij}(p_i,p_j)\|}^2_{W_{ij}},
\end{equation*}
where $E_{\text{odom}}$ represents all pairs $(i,j)$ from the consecutive measurements, $E_{\text{loop}}$ all pairs from the loop closure measurements, and $\|\cdot\|_{W_{ij}}$ denotes the norm induced by the weight matrix $W_{ij}$.

From~\eqref{eq:slam-relative-pose}, it follows that the Jacobian of $\Phi(p_{i}, p_{j})$ is given by
\begin{equation*}
  \begin{aligned}
    & D\Phi(p_{i}, p_{j}) \\[0.5ex]
    = &
      \begin{pmatrix}
        \frac{\partial\Phi_{1}}{\partial x_{i}} & \frac{\partial\Phi_{1}}{\partial y_{i}} & \frac{\partial\Phi_{1}}{\partial \theta_{i}} & \frac{\partial\Phi_{1}}{\partial x_{j}} & \frac{\partial\Phi_{1}}{\partial y_{j}} & \frac{\partial\Phi_{1}}{\partial \theta_{j}} \\[0.5ex]
        \frac{\partial\Phi_{2}}{\partial x_{i}} & \frac{\partial\Phi_{2}}{\partial y_{i}} & \frac{\partial\Phi_{2}}{\partial \theta_{i}} & \frac{\partial\Phi_{2}}{\partial x_{j}} & \frac{\partial\Phi_{2}}{\partial y_{j}} & \frac{\partial\Phi_{2}}{\partial \theta_{j}} \\[0.5ex]
        \frac{\partial\Phi_{3}}{\partial x_{i}} & \frac{\partial\Phi_{3}}{\partial y_{i}} & \frac{\partial\Phi_{3}}{\partial \theta_{i}} & \frac{\partial\Phi_{3}}{\partial x_{j}} & \frac{\partial\Phi_{3}}{\partial y_{j}} & \frac{\partial\Phi_{3}}{\partial \theta_{j}}
      \end{pmatrix}
    \\[0.5ex]
    = &
        \setlength\arraycolsep{4pt}
    \begin{pmatrix}
      -\cos\theta_{i} & -\sin\theta_{i} &  -\sin\theta_{i}(x_j-x_i) + \cos\theta_{i}(y_{j}-y_{i}) & \cos\theta_{i}  & \sin\theta_{i} & 0 \\[0.5ex]
      \sin\theta_{i}  & -\cos\theta_{i} &  -\cos\theta_{i}(x_j-x_i) - \sin\theta_{i}(y_{j}-y_{i}) & -\sin\theta_{i} & \cos\theta_{i} & 0 \\[0.5ex]
          0          &       0         & -1 &     0           &      0        & 1
    \end{pmatrix},
  \end{aligned}
\end{equation*}
where we dropped the arguments in the second row for simplicity.  The gradient of $\mathcal{J}$ can be written as
\begin{equation*}
  \nabla\mathcal{J}(p_{1},\ldots,p_{N}) =
  \begin{pmatrix}
    \nabla_{p_{1}}\mathcal{J}(p_{1},\ldots,p_{N}) \\[0.5ex] \vdots \\[0.5ex] \nabla_{p_{N}}\mathcal{J}(p_{1},\ldots,p_{N})
  \end{pmatrix}\in\R^{3N},
\end{equation*}
where
\begin{equation}
  \label{eq:slam-gradient}
  \begin{aligned}
    \nabla_{p_{i}}\mathcal{J}(p_{1},\ldots,p_{N})
    &=
    \begin{pmatrix}
      \frac{\partial\mathcal{J}}{\partial x_{i}}(p_{1},\ldots,p_{N}) \\[0.5ex]
      \frac{\partial\mathcal{J}}{\partial y_{i}}(p_{1},\ldots,p_{N}) \\[0.5ex]
      \frac{\partial\mathcal{J}}{\partial \theta_{i}}(p_{1},\ldots,p_{N})
    \end{pmatrix}
    \\[0.5ex]
  &= \sum_{\substack{(k,\ell)\in E,\\\color{black} \{k,\ell\}\ni i}} D_{\color{black} p_i}\Phi(p_{k}, p_{\ell})^{T}\,W_{k\ell}\,r_{k\ell}(p_{k},p_{\ell}).
  \end{aligned}
\end{equation}
The Hessian $H$ of $\mathcal{J}$ can be written as
\begin{equation*}
  H =
  \begin{pmatrix}
    H_{p_{1},p_{1}} & \ldots & H_{p_{1},p_{N}} \\[0.5ex]
    \vdots        &        & \vdots        \\[0.5ex]
    H_{p_{N},p_{1}} & \ldots & H_{p_{N},p_{N}}
  \end{pmatrix},
\end{equation*}
where $H_{p_{i},p_{j}}$ is given by
\begin{equation}
  \label{eq:slam-hessian}
  \begin{aligned}
    H_{p_{i},p_{j}} =\ & \nabla_{p_{i},p_{j}}^{2}\mathcal{J}(p_{1},\ldots,p_{N}) \\[0.5ex]
    =\ & \underbrace{\sum_{\substack{(k,\ell)\in E,\\
    \color{black} \{k,\ell\}\ni i,j}}
    D_{\color{black}p_i}
    \Phi(p_{k}, p_{\ell})^{T}W_{k\ell}D_{\color{black}p_\ell}\Phi(p_{k}, p_{\ell})}_{\text{Gauss--Newton approximation}} + B_{p_{i},p_{j}},
  \end{aligned}
\end{equation}
where the matrix $B_{p_{i},p_{j}}$ contains the remaining second order terms.
These terms are neglected in Gauss-Newton.

From~\eqref{eq:slam-hessian}, it follows that
\begin{equation*}
  H = H_{\text{GN}} + B,
\end{equation*}
where $H_{\text{GN}}$ refers to the Gauss--Newton approximation of $H$. Here, $B$ refers to the second order derivatives of $\Phi$, which are neglected in Gauss--Newton.

To derive a more compact matrix notation often used in the literature, let
\[
r(p) =
\begin{pmatrix}
r_{i_1 j_1}(p_{i_1},p_{j_1}) \\
\vdots \\
r_{i_m j_m}(p_{i_m},p_{j_m})
\end{pmatrix}
\in \mathbb{R}^{3m}
\]
denote the stacked residual vector containing all edge residuals, where
$m = |E|$ is the number of edges in the pose graph.
Furthermore, let
\[
W = \operatorname{diag}(W_{i_1 j_1},\ldots,W_{i_m j_m})
\]
be the block diagonal matrix collecting all weight matrices.
Then the objective function~\eqref{eq:slam-obj} can be written as
\[
\mathcal J(p) = \frac12\, r(p)^T W r(p).
\]

Let
\[
J(p) = \frac{\partial r(p)}{\partial p}
\]
denote the Jacobian of the stacked residual vector.
Then the gradient of $\mathcal J$ can be written as
\[
\nabla \mathcal J(p) = J(p)^T W r(p).
\]

Finally, the Gauss--Newton approximation of the Hessian is given by
\begin{equation}
H_{\mathrm{GN}}(p) = J(p)^T W J(p).
\label{eq:hgn}
\end{equation}

For symmetric positive definite  weight matrices $W_{ij}$, and after
the natural null space of the problem is removed (e.g., by fixing one pose),
$H_{\mathrm{GN}}$ is symmetric positive definite, and the conjugate gradient method is applicable.

Throughout this paper, we will use Gauss--Newton, and we will write
\begin{equation}
A:=H_{GN}
\label{eq:A}
\end{equation}


\subsection{Local-to-Global Assembly in Pose Graph SLAM}
The algebraic structure of the Gauss--Newton matrix corresponds
to the assembly process known from finite element discretizations.  Each
edge \((i,j)\in E\) of the pose graph defines a residual depending only
on the two incident poses.  After linearization, this residual contributes
a small local matrix
\[
  H_{ij} = J_{ij}^T W_{ij} J_{ij}
\]
to the global Gauss--Newton matrix.  If \(R_{ij}\) denotes the restriction
to the degrees of freedom associated with poses \(i\) and \(j\), then
the global matrix can be written in the assembled form
\[
  H =
  \sum_{(i,j)\in E}
  R_{ij}^T H_{ij} R_{ij}.
\]
This is analogous to finite element assembly, where local element
stiffness matrices are accumulated into a global sparse stiffness matrix.
The analogy is useful for domain decomposition methods: subdomains can
be defined by subsets of graph edges or vertices, and local subdomain
matrices are obtained by restricting the global problem to these
subgraphs.

\section{Overlapping Schwarz Domain Decomposition Preconditioners}

The general concept of domain decomposition methods (DDM) is based on the divide-and-conquer paradigm: A large problem is decomposed into multiple smaller, more manageable subproblems which can be solved separately and in parallel. A comprehensive introduction to DDM can be found in~\cite{Toselli:2005:DDM,sbg_book_1996}. DDM are typically accelerated by Krylov methods such as conjugate gradients.

The divide-and-conquer approach makes DDM very attractive for the parallel solution of large, ill-conditioned problems. This motivated significant research activity on the usage of DDM in large  problems from computational fluid mechanics~\cite{gropp1992domain}, structural mechanics~\cite{gosselet2006non,itersub_2009,heinlein:2023:PST}, or coupled problems~\cite{corigliano2013domain,Kiefer:2023:MPO,Klawonn:2020:CHM}.

Although DDM were originally developed for the solution of large systems from the discretization of partial differential equations, e.g., by finite elements, the underlying ideas can also be formulated algebraically.

%



Overlapping Schwarz methods are a subclass within DDM, which use overlapping subdomains. For additive one-level overlapping Schwarz, local corrections are computed on each subdomain and then combined additively to a global approximate solution; see~\eqref{eq:OS-preconditioner}. Two-level methods add a mechanism for fast global transport of information.

{From a parallel computing perspective, each subdomain typically corresponds
to a process or a group of processes as in~\cite{heinlein:2023:PST}, or more
generally to a thread or processing unit in a parallel implementation.
The restriction operators $R_i$ and their transposes then represent parallel
communication or data movement patterns, e.g., in message-passing (MPI) or
shared-memory environments. Applying $R_i$ corresponds to gathering local
data, while applying $R_i^T$ corresponds to scattering local corrections
back to the global vector. The overlap determines the amount of communication
or synchronization required.
Such data movement patterns occur naturally not only on distributed-memory systems, but also on shared-memory architectures and embedded multi-core platforms commonly used in robotics.

}

\subsection{Graph-Based Domain Decomposition}
After Gauss--Newton linearization, pose-graph SLAM leads to sparse symmetric positive definite linear systems
\begin{equation}
A x = b,
\end{equation}
where
$ 
A = J^T W J,
$ 
where \(J\) denotes the Jacobian of the residuals and \(W\) is the weighting matrix; see Section~\ref{seq:lsq}.

In this work, we will consider a simple model problem in which a robot repeatedly traverses a unit square as a benchmark problem in Section~\ref{sec:modelproblem}. Each traversal generates one loop in the pose graph. The resulting graph therefore consists of multiple strongly-connected trajectory segments coupled through loop closures; see also Fig.~\ref{fig:optimized_results}.


As the number of loops increases, the linear systems become increasingly difficult to solve iteratively. We therefore employ overlapping Schwarz domain decomposition preconditioners in combination with conjugate gradients.

Since pose-graph SLAM naturally leads to sparse graph structures rather than geometric meshes, we formulate the Schwarz methods directly on the graph.
%
%

\begin{figure}[tb]
\begin{center}
\begin{tikzpicture}[line join=round]
  \basechain

  \foreach \i in {1,...,5} {\coordinate (a\i) at ({\i-1},1.10);}
  \foreach \i in {5,...,9} {\coordinate (b\i) at ({\i-1},1.85);}

  \draw[subedgeA] (a1) -- (a2) -- (a3) -- (a4) -- (a5);
  \draw[subedgeB] (b5) -- (b6) -- (b7) -- (b8) -- (b9);

  \foreach \i in {1,...,5} {\node[nodeA] at (a\i) {};}
  \foreach \i in {5,...,9} {\node[nodeB] at (b\i) {};}

  \node[ringA] at (a5) {};
  \node[ringB] at (b5) {};
  \draw[conn] (p5) -- (a5);
  \draw[conn] (p5) -- (b5);

  \node[ilabel, SegA, above=5pt] at ($(a1)!0.5!(a5)$) {$\mathcal{I}_{1}$};
  \node[ilabel, SegB, above=5pt] at ($(b5)!0.5!(b9)$) {$\mathcal{I}_{2}$};

\end{tikzpicture}
\caption{Nonoverlapping edge partition. The node index sets meet in the interface pose $p_5$.
\label{fig:nonoverlapping}
}
\end{center}
\end{figure}

\begin{figure}[tb]
\begin{center}
\begin{tikzpicture}[line join=round]
  \basechain

  \foreach \i in {1,...,6} {\coordinate (a\i) at ({\i-1},1.10);}
  \foreach \i in {4,...,9} {\coordinate (b\i) at ({\i-1},1.85);}

  \draw[subedgeA] (a1) -- (a2) -- (a3) -- (a4) -- (a5) -- (a6);
  \draw[subedgeB] (b4) -- (b5) -- (b6) -- (b7) -- (b8) -- (b9);

  \foreach \i in {1,...,6} {\node[nodeA] at (a\i) {};}
  \foreach \i in {4,...,9} {\node[nodeB] at (b\i) {};}

  \foreach \i in {4,5,6} {
    \node[ringA] at (a\i) {};
    \node[ringB] at (b\i) {};
    \draw[conn] (p\i) -- (a\i);
    \draw[conn] (p\i) -- (b\i);
  }

  \node[ilabel, SegA, above=5pt] at ($(a1)!0.5!(a6)$) {$\mathcal{I}_{1}^{\delta}$};
  \node[ilabel, SegB, above=5pt] at ($(b4)!0.5!(b9)$) {$\mathcal{I}_{2}^{\delta}$};

\end{tikzpicture}
\caption{Construction of the overlapping subdomains in the case of overlap one: each subdomain is enlarged by one neighboring graph edge; the common nodes are now $p_4,p_5,p_6.$ The interface node is $p_5$.
\label{fig:overlapone}
}
\end{center}
\end{figure}

We first define a nonoverlapping edge decomposition
by partitioning the edge set into nonoverlapping subdomains
\begin{equation}
E
=
\dot\bigcup_{i=1}^N E_i.
\end{equation}
In our model problem, each robot loop defines one edge subdomain \(E_i\).

The decomposition is therefore nonoverlapping with respect to measurement edges: every edge belongs to exactly one subdomain.

For every edge subdomain \(E_i\), we define the associated vertex set
\begin{equation}
V_i
=
\left\{
v\in V :
v \text{ is incident to an edge in } E_i
\right\}.
\end{equation}

Although the edge sets \(E_i\) are disjoint,
neighboring subdomains will share an interface pose, i.e.,
\[
V_i \cap V_j \neq \emptyset;
\]
see Fig.~\ref{fig:nonoverlapping}.
This is analogous to finite element domain decomposition, where subdomains are nonoverlapping in their interiors but their closures meet at the common interface.



The algebraic degrees of freedom associated with subdomain \(i\) are given by the index set
\begin{equation}
\mathcal I_i
=
\left\{
\text{degrees of freedom associated with vertices in } V_i
\right\}.
\end{equation}



For the construction of the overlap, we distinguish between odometry edges and loop-closure edges,
$E = E_{\mathrm{odom}} \cup E_{\mathrm{loop}}$.
The overlapping subdomains are obtained by enlarging each subdomain by one layer of neighboring graph edges. In this enlargement step, loop-closure edges are not followed as graph-neighborhood edges. However, once the overlapping vertex sets have been constructed, all matrix couplings between degrees of freedom contained in the corresponding overlapping index set are included in the local matrices. Thus, loop-closure couplings whose endpoints lie in the overlap are retained in the local Schwarz problems.

We always use an overlap of one, i.e., each subdomain is enlarged by one layer of neighboring graph edges; see also Fig.~\ref{fig:overlapone}.

More precisely, we define
\begin{equation}
E_i^\delta
=
E_i
\cup
\left\{
e\in E_{\mathrm{odom}} :
e \text{ is incident to a vertex in } V_i
\right\}.
\end{equation}
The enlarged vertex set becomes
\begin{equation}
V_i^\delta
=
\left\{
v\in V :
v \text{ is incident to an edge in } E_i^\delta
\right\}.
\end{equation}
The corresponding overlapping algebraic index set is denoted by
$
\mathcal I_i^\delta.
$


\subsection{One-Level Additive Schwarz Method}

For each overlapping subdomain defined by the index set \(\mathcal I_i^\delta\), we define the restriction operator
\[
R_i : \mathbb R^n \to \mathbb R^{n_i},
\]
which extracts the unknowns belonging to the corresponding index set; cf.~\cite{Toselli:2005:DDM,sbg_book_1996}.
The corresponding local Schwarz matrix is given by
\begin{equation}
A_i = R_i A R_i^T.
\end{equation}
The one-level additive Schwarz preconditioner~\cite{Toselli:2005:DDM,sbg_book_1996}
is defined as
\begin{equation}
M^{-1}_{AS}
=
\sum_{i=1}^N
R_i^T A_i^{-1} R_i.
\label{eq:OS-preconditioner}
\end{equation}

Algorithmically, in each CG iteration, an application of the one-level preconditioner consists of:
\begin{enumerate}[label=\arabic*.]
\item Restricting the residual to all subdomains,
\item Solving independent local problems,
\item Adding the local corrections.
\end{enumerate}
Since all local solves are independent, the method is naturally parallel.


One-level Schwarz methods propagate information only through neighboring subdomains.
As the number of robot loops increases, information about global trajectory corrections must travel successively from one subdomain to the next. Consequently, the number of CG iterations grows with the number of subdomains.
This phenomenon is well known for elliptic PDEs and also appears in SLAM model problems. To obtain scalable convergence, an additional global coarse correction is required.

\subsection{Two-Level Additive Schwarz Method}

The two-level additive Schwarz preconditioner~\cite{Toselli:2005:DDM,sbg_book_1996}
augments the local corrections by a global coarse correction
for fast global transport of information,
\begin{equation}
M^{-1}_{TL}
=
\sum_{i=1}^N
R_i^T A_i^{-1} R_i
+
R_0^T A_0^{-1} R_0.
\end{equation}
Here,
\begin{equation}
A_0 = R_0 A R_0^T
\end{equation}
denotes the coarse problem.
The coarse correction enables global information to propagate through the entire pose graph within a single iteration. Classically the coarse correction is constructed using a coarse representation of the original problem, e.g., using a coarse grid~\cite{Toselli:2005:DDM,sbg_book_1996}.

\subsection{Algebraic Construction of a GDSW-type Coarse Basis}
An advantage of the GDSW (Generalized Dryja--Smith--Widlund) coarse basis functions~\cite{Dohrmann:2008:FEM,Dohrmann:2008:DDL,Heinlein:2016:PIT,heinlein:2023:PST} for overlapping Schwarz methods is that numerical and parallel scalability
can be achieved~\cite{heinlein:2017:POS,heinlein:2023:PST} using only information extracted from the system matrix, or, in some cases, using little additional geometric information, e.g., to construct rotations.

In GDSW methods for elliptic partial differential equations in two dimensions, the interface functions are defined using the nullspace $Z$ of the global Neumann matrix to the vertices and edges, which form a nonoverlapping decomposition of the interface $\Gamma$ of the nonoverlapping domain decomposition.

In our case, although the problem is two-dimensional, the interface of neighboring nonoverlapping subdomains consists of exactly one robot pose.
Therefore, every interface consists of a single pose with the three unknowns
$
(x,y,\phi).
$
Following the GDSW approach, for every interface pose, we define the three columns of the interface function $\Phi_\Gamma$
\begin{equation}
\phi^T_{\Gamma,1}=(1,0,0)^T,
\qquad
\phi^T_{\Gamma,2}=(0,1,0)^T,
\qquad
\phi^T_{\Gamma,3}=(0,0,1)^T.
\end{equation}


In two-dimensional mechanics, GDSW typically makes use of
linearized rotations $(-y,x,1)$ in the coarse space to obtain numerical and parallel scalability, e.g.,~\cite{Heinlein:2016:PIT,heinlein:2017:POS,heinlein:2023:PST}. However, this choice would not be algebraic, i.e., $(-y,x,1)$ cannot be constructed from the system matrix, only.
Therefore, our choice is $(0,0,1)^T$, which is algebraic and naturally arises if the centers of all rotations are moved to the corresponding interface node.


The interface basis vectors are extended into the neighboring subdomains by discrete harmonic extension.
Let
$
\phi_\Gamma
$
denote prescribed interface values.
After partitioning the local matrix into interior and interface degrees of freedom,
\begin{equation}
A_i
=
\begin{pmatrix}
A_{II} & A_{I\Gamma} \\
A_{\Gamma I} & A_{\Gamma\Gamma}
\end{pmatrix},
\end{equation}
the interior values $\phi_I$ are computed from 
\begin{equation}
A_{II} \phi_I
=
-
A_{I\Gamma} \phi_\Gamma.
\end{equation}

\begin{figure}[t!]
\begin{minipage}[t]{1.0\columnwidth}
    \centering
    \includegraphics[width=\linewidth]{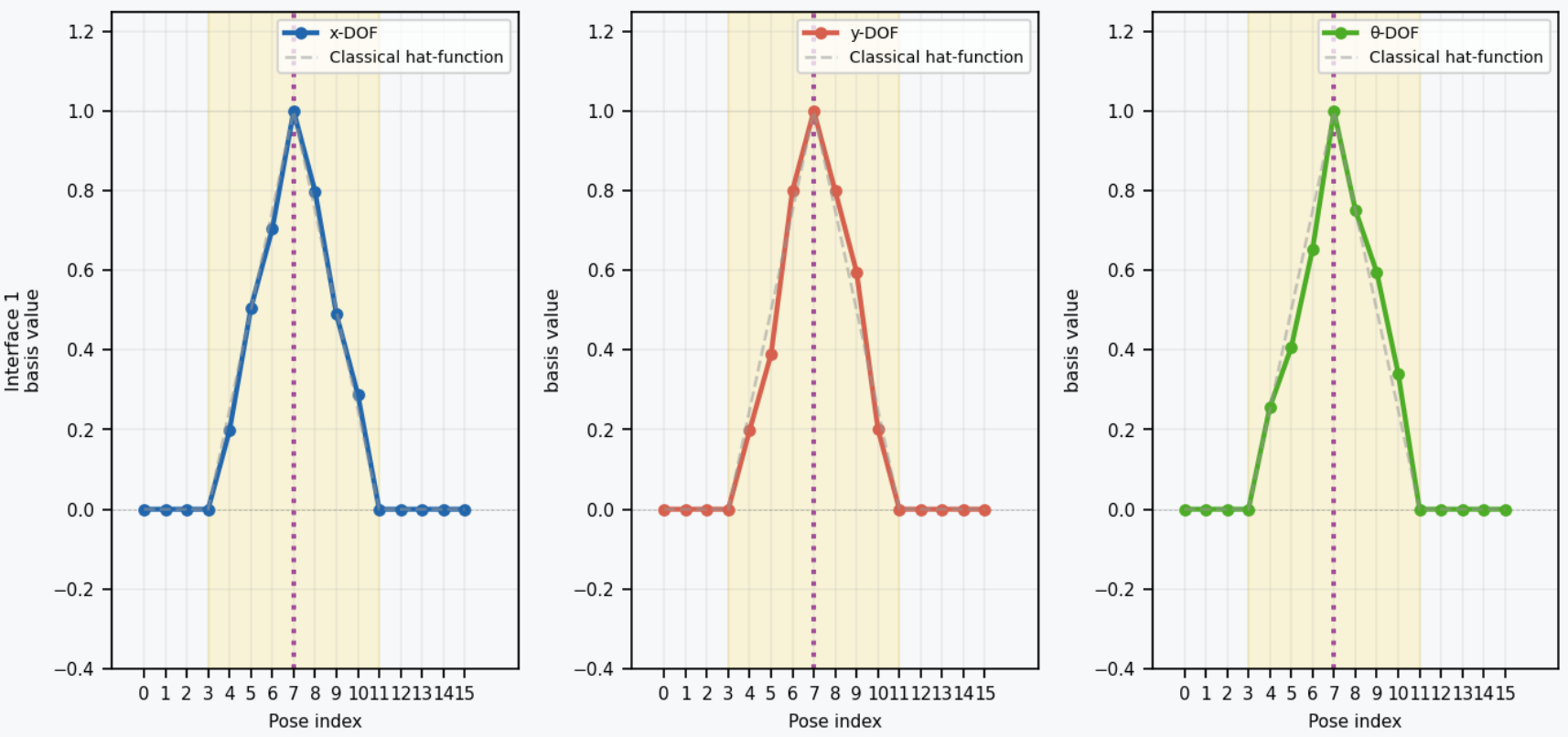}
    \caption{An example of our GDSW-type coarse functions for the benchmark SLAM problem described in Section~\ref{sec:numer}. 
    The classical hat-functions are overlaid as comparison.
    Note that in one-dimension, the GDSW functions for the Laplace equation are the classical hat-functions.
    }\label{fig:hat}
\end{minipage}\hfill
\end{figure}

The resulting coarse basis functions are energy-minimal 
extensions of the interface modes into the neighboring subdomains; see Figure~\ref{fig:hat}.
Note that the coarse basis functions are quite close to standard hat-functions but not identical in our SLAM benchmark problem, since the trajectory contains both collinear odometry segments and right-angle turns.


We collect all $\phi_\Gamma$ in the coarse basis matrix
$
\Phi.
$
Using the coarse restriction operator
$
R_0 = \Phi^T,
$
the coarse matrix becomes
\begin{equation}
A_0
=
R_0 A R_0^T
=
\Phi^T A \Phi.
\end{equation}

The resulting coarse space captures the dominant global error components responsible for slow convergence of the one-level method.

\subsection{Interpretation for Pose-Graph SLAM}

In the SLAM-context, the two-level Schwarz method has a natural interpretation:
\begin{itemize}
\item local Schwarz corrections optimize local trajectory consistency,
\item the coarse space propagates global loop-closure information through the entire graph.
\end{itemize}
This significantly improves convergence as the number of loops increases.
Our numerical experiments demonstrate that
\begin{itemize}
\item unpreconditioned CG scales poorly,
\item one-level Schwarz improves convergence but is not scalable,
\item the GDSW coarse space yields bounded iteration counts for increasing number of subdomains.
\end{itemize}


\section{Numerical Results and Discussion}
\label{sec:modelproblem}
To investigate the convergence of domain decomposition methods in SLAM, we consider the optimization of a pose graph SLAM problem in two dimensions. 

Artificial odometry data are generated by adding Gaussian noise to ground truth data, and optimization is subsequently performed on this noisy data.
%
In our test case, the robot moves from the origin following a square trajectory, moving one unit to the right and turning left by 90\textdegree, repeating the motion until it returns to the starting point. We denote this movement around the unit square as \textit{one loop}.


When the robot revisits the starting position, it recognizes that it has returned to a previously visited location, and a loop closure edge is added to the pose graph; see Fig.~\ref{fig:slam-loop-closure}. In our experiments, the weight matrix is $W_{ij}=100I_3$ if the edge $(i,j)$ corresponds to a loop closure, and $W_{ij}=20I_3$ otherwise. The higher weight reflects the higher confidence assigned to loop closure constraints compared to the other measurements.

We solve the least-squares problem using Gauss--Newton, where the sparse  systems from linearization are solved by conjugate gradients using additive overlapping Schwarz as a preconditioner.
A relative tolerance of $10^{-8}$ for conjugate gradients is used for the solution of the sparse systems. As for the least-squares problem, the stopping criterion is twofold: an absolute tolerance of $10^{-8}$, or a relative decrease of the gradient by $10^{-6}$.
Since the inner conjugate gradient solver uses a relative stopping criterion, it is less sensitive to the problem size. However, the absolute stopping criterion in the outer nonlinear iteration may still lead to problem-size-dependent behavior.

In the subsequent text and tables, we will write CG for conjugate gradients and GN for Gauss--Newton.

\begin{figure}[t!]
\begin{minipage}[t]{0.4\columnwidth}
    \centering
    \includegraphics[width=\linewidth]{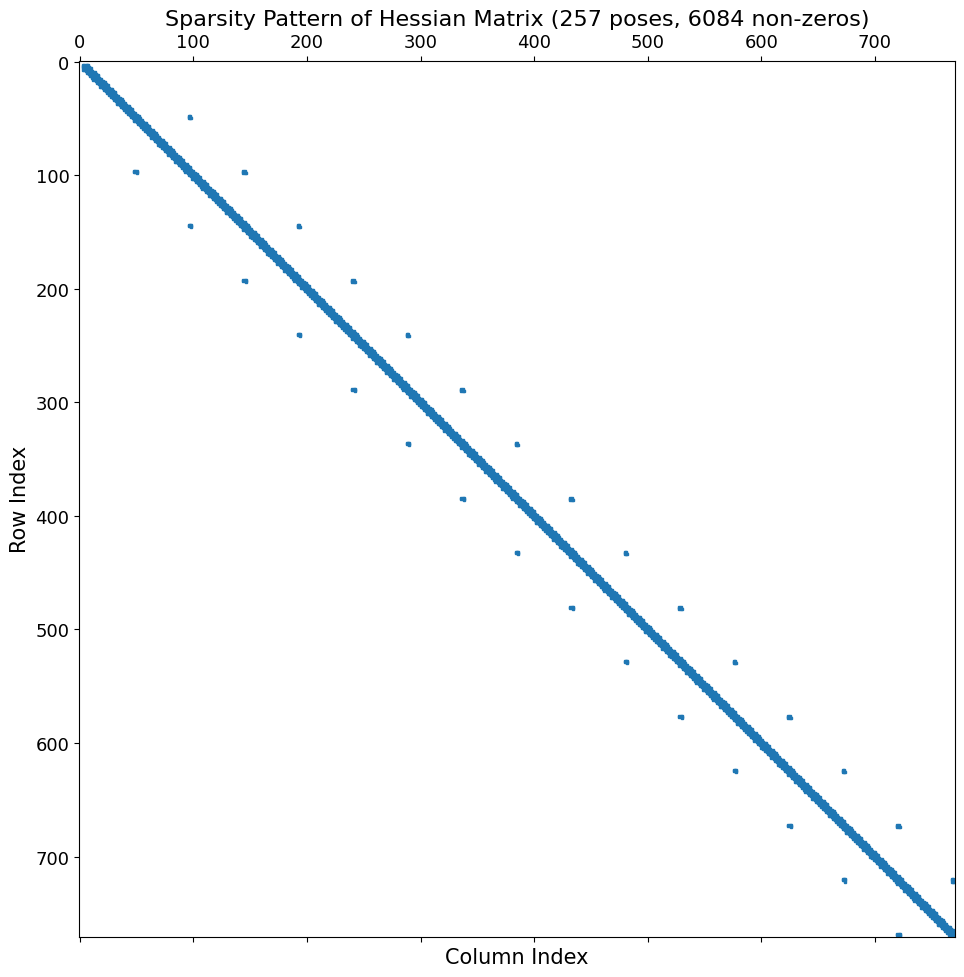}
    \caption{Sparsity pattern of the system matrix. Loop closures introduce additional off-diagonal coupling visible in the sparsity.}\label{fig:sparsity}
\end{minipage}\hfill
\begin{minipage}[t]{0.5\columnwidth}
    \centering
    \includegraphics[width=\linewidth]{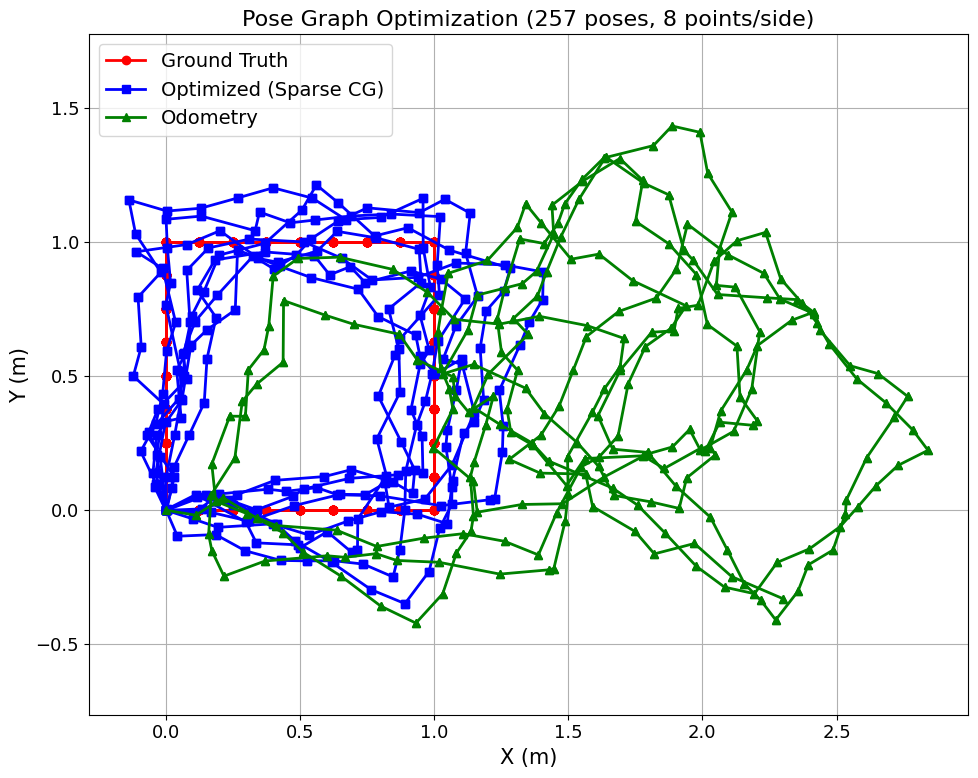}
    \caption{Comparison of an optimized robot path with odometry and ground truth; 8 loops, 8 points per side.}\label{fig:optimized_results}
\end{minipage}
\end{figure}

\label{sec:numer}
In our experiments, each loop defines a subdomain for overlapping Schwarz. We use minimal overlap, where the overlap consists of three degrees of freedom. In our choice of the subdomains, the overlap includes the off-diagonal coupling introduced by the loop closures; see Fig.~\ref{fig:sparsity}.

In other words, the degrees of freedom connected by the loop-closure constraints are contained in the overlapping index sets in this test configuration, so that the corresponding matrix couplings are included in the local Schwarz problems.


To evaluate the robustness and numerical scalability of this overlapping Schwarz preconditioner, the problem size is 
increased by increasing 
the number of loops (and thus the number of subdomains) traversed by the robot. We also consider the case where the number of odometry points on each side of the unit square is increased.

An example of the sparsity pattern of our system matrix
is shown in Fig.~\ref{fig:sparsity}.
Here, 8 loops are traversed, with 8 points on each side of the unit square. Fig.~\ref{fig:sparsity} illustrates the sparsity of $H$ and highlights the off-diagonal blocks resulting from the long range coupling introduced by the loop closures.

Fig.~\ref{fig:optimized_results} shows the solution of the nonlinear least-squares problem for the case where 8 cycles of unit squares have been traversed by the robot, with 8 odometry points taken at each side of the unit square. The solution of the least-squares problem converges to the same results with and without the use of our preconditioner.

\begin{table}[tb!] 
\centering
\resizebox{0.8\textwidth}{!}
{ \begin{tabular}{cc||cc|cc|cc} 
\toprule \multicolumn{2}{c||}{} & \multicolumn{2}{c|}{\textbf{None}} & \multicolumn{2}{c}{\textbf{One-Level}} & \multicolumn{2}{c}{\textbf{GDSW}} \\
\cmidrule(lr){3-4}\cmidrule(lr){5-6}\cmidrule(lr){7-8} \textbf{Loops} & \textbf{Points/Side} & \textbf{GN} & \textbf{CG} & \textbf{GN} & \textbf{CG} & \textbf{GN} & \textbf{CG} \\
\midrule
4 & 4 & 5 & 153.0 & 5 & 16.0 & 5 & 13.0 \\
& 8 & 5 & 273.6 & 5 & 15.8 & 5 & 12.4 \\
& 16 & 6 & 510.8 & 6 & 16.0 & 6 & 12.3 \\
& 32 & 5 & 1036.2 & 5 & 15.8 & 5 & 12.4 \\
& 64 & 6 & 2115.7 & 6 & 15.8 & 6 & 12.5 \\
& 128 & 8 & 3750.0 & 8 & 14.6 & 8 & 11.6 \\
\midrule
8 & 4 & 5 & 206.0 & 5 & 25.2 & 5 & 14.8 \\
& 8 & 6 & 385.3 & 6 & 24.7 & 6 & 14.5 \\
& 16 & 6 & 705.5 & 6 & 24.0 & 6 & 14.5 \\
& 32 & 5 & 1342.2 & 5 & 24.6 & 5 & 14.6 \\
& 64 & 7 & 3000.4 & 7 & 23.9 & 7 & 14.4 \\
& 128 & 8 & 5441.2 & 8 & 22.0 & 8 & 13.4 \\
\midrule
16 & 4 & 6 & 329.3 & 6 & 41.5 & 6 & 15.8 \\
& 8 & 5 & 551.8 & 5 & 40.6 & 5 & 16.0 \\
& 16 & 6 & 993.2 & 6 & 40.0 & 6 & 15.3 \\
& 32 & 6 & 2031.7 & 6 & 37.2 & 6 & 15.3 \\
& 64 & 7 & 3694.0 & 7 & 38.1 & 7 & 15.3 \\
& 128 & 9 & 9158.8 & 9 & 40.7 & 9 & 15.1 \\
\midrule
32 & 4 & 6 & 549.2 & 6 & 73.2 & 6 & 16.7 \\
& 8 & 6 & 898.7 & 6 & 73.0 & 6 & 16.7 \\
& 16 & 7 & 1557.3 & 7 & 65.7 & 7 & 16.7 \\
& 32 & 6 & 2567.3 & 6 & 66.8 & 6 & 16.2 \\
& 64 & 8 & 8014.6 & 8 & 75.0 & 8 & 16.2 \\
& 128 &16 &12549.7 &16 & 61.8 &16 & 15.1 \\
\bottomrule
\end{tabular} }
\caption{
Comparison of the number of
conjugate gradient (CG) iterations for the cases
without using a preconditioner, when using the one-level additive overlapping Schwarz preconditioner, and when using the two-level preconditioner.
We vary the number of loops and the points per side. The number of loops is identical to the number of overlapping subdomains, and the subdomain size increases with the points per side.
}\label{tab:numerical_performance}
\end{table}

\begin{table}[tb!]
\centering
\resizebox{0.8\textwidth}{!}
{ \begin{tabular}{cc||cc|cc|cc}
\toprule
\multicolumn{2}{c||}{} & \multicolumn{2}{c|}{\textbf{None}} & \multicolumn{2}{c}{\textbf{One-Level}} & \multicolumn{2}{c}{\textbf{GDSW}} \\
\cmidrule(lr){3-4}\cmidrule(lr){5-6}\cmidrule(lr){7-8} \textbf{Loops} & \textbf{Points/Side} & \textbf{GN} & \textbf{CG} & \textbf{GN} & \textbf{CG} & \textbf{GN} & \textbf{CG} \\
\midrule
4 & 16 & 6 & 510.8 & 6 & 16.0 & 6 & 12.3 \\
8 & 16 & 6 & 705.5 & 6 & 24.0 & 6 & 14.5 \\
16 & 16 & 6 & 993.2 & 6 & 40.0 & 6 & 15.3 \\
32 & 16 & 7 &1557.3 & 7 & 65.7 & 7 & 16.7 \\
64 & 16 & 6 &2370.7 & 6 &121.5 & 6 & 16.7 \\
128 & 16 & 9 &6633.2 & 9 &264.3 & 9 & 16.8 \\
\bottomrule
\end{tabular} }
\caption{Numerical scalability for an increasing number of subdomains. The averaged number of additive Schwarz preconditioner CG iterations remains bounded for an increasing number of subdomains when a two-level preconditioner is used}\label{tab:weak_numerical}
\end{table}

The numerical results obtained from the tests performed are presented in Table~\ref{tab:numerical_performance}. In our numerical experiments, the robot traverses 4, 8, 16, and 32 loops. The number of loops coincides with the number of overlapping Schwarz subdomains, which means that at most 32 subdomains are used. Thus, this tests numerical scalability with respect to the number of subdomains (weak scalability). For each of these cases, we increase the number of odometry points at each side of the unit square until 128 points per side. This tests the numerical scalability for the case of increasing subdomain size. Note that in all cases minimal overlap is used.

In Table~\ref{tab:numerical_performance}, we can see that the preconditioner drastically reduces the number of conjugate gradient (CG) iterations needed to reach convergence: while the number of unpreconditioned CG iterations increases to more than $10\,000$ Krylov iterations, this number is greatly reduced for the one-level preconditioner case, even though an increase in iterations typical of the one-level preconditioner could still be seen with increased number of subdomains. The same observation is not seen in the two-level preconditioner case, where the increase in iterations is much less noticeable.

In addition, the number of CG iterations does not increase significantly with the number of subdomains or with the subdomain size for the two-level case. This shows good numerical scalability for this particular test case, i.e., despite the increase in problem size, the number of CG iterations remains essentially constant, and the condition numbers
$\lambda_{\rm max}/\lambda_{\rm min}$ of the preconditioned system stay bounded.

\begin{figure}[t]
  \centering
  \includegraphics[width=1.0\linewidth]{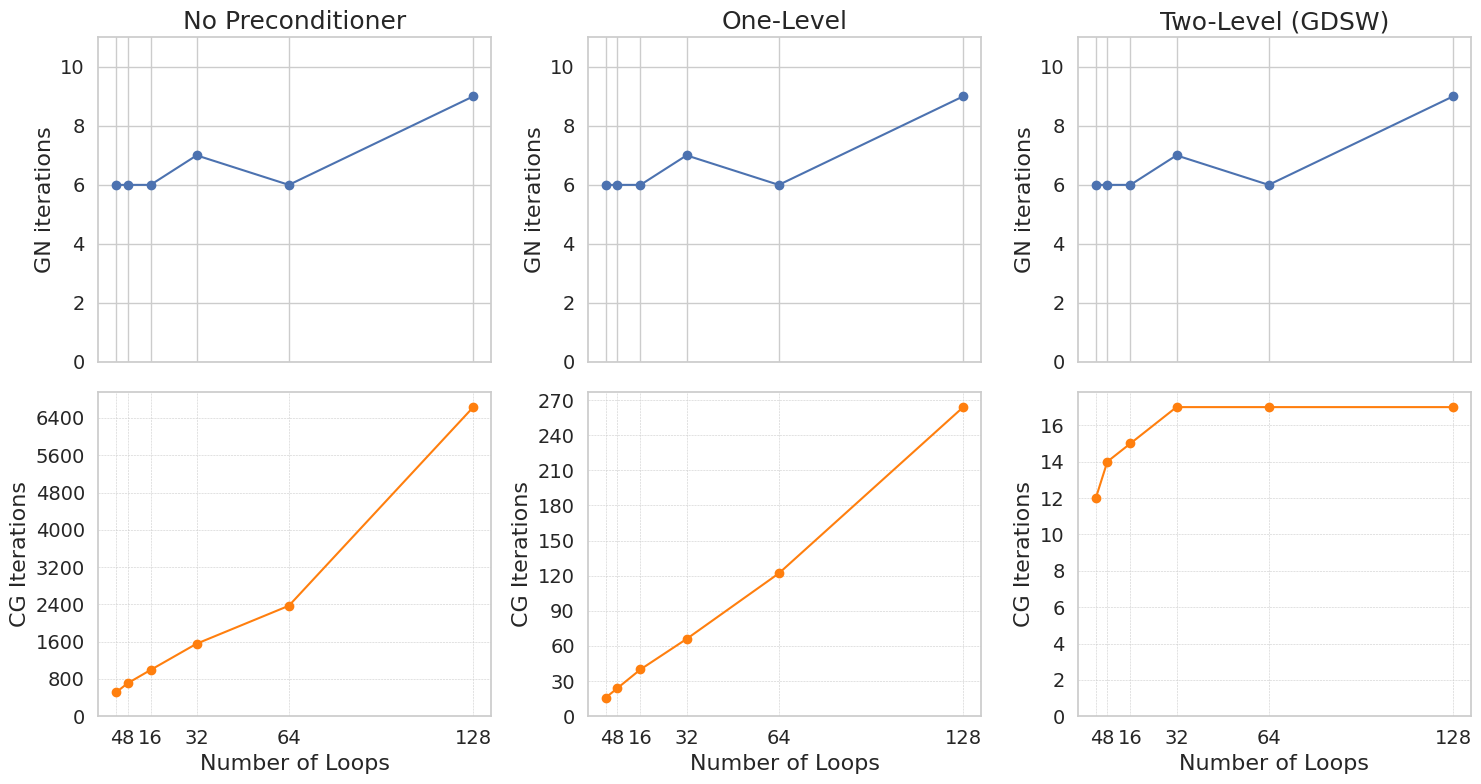}
  \caption{Comparison of numerical scalability with and without using the overlapping Schwarz preconditioners by varying the number of loops traversed by the robot while keeping the number of points at each side of the unit square constant at 16. The number of subdomains is always chosen equal to the number of loops.
    \\
  }\label{fig:scalability}
\end{figure}

The numerical scalability with respect to the number of subdomains (and loops) is presented in Table~\ref{tab:weak_numerical} and Fig.~\ref{fig:scalability}. We observe, again, that the number of CG iterations remains bounded, when the two-level additive Schwarz preconditioner is used. This is in contrast to the case without preconditioner.



Let us briefly discuss the number of nonlinear iterations. The numerical results in Table~\ref{tab:numerical_performance} indicate that the number of Gauss--Newton iterations somewhat increases with the number of points per loop. To some extent, this may be a consequence of the absolute stopping criterion for Gauss--Newton, which is based on the standard Euclidean norm.
However, for a fixed number of points per loop, the number of Gauss--Newton iterations appears to remain bounded as the number of loops increases. In Table~\ref{tab:weak_numerical}, however, we observe an increase from 6 to 9 Gauss--Newton iterations when scaling to 128 loops. Further studies are necessary to determine whether the convergence of Gauss--Newton deteriorates for this problem as the number of loops grows.

To assess the importance of the rotational coarse mode, Table~\ref{tab:weak_null} compares the GDSW-type coarse space with and without the rotational component in the local null space.

\begin{table}[tb!]
\centering
\resizebox{1.0\textwidth}{!}
{ \begin{tabular}{cc||cc|cc}
\toprule
\multicolumn{2}{c||}{}  & \multicolumn{2}{c|}{\textbf{GDSW without Rotation}} & \multicolumn{2}{c}{\textbf{GDSW with Rotation}} \\
\cmidrule(lr){3-4}\cmidrule(lr){5-6} \textbf{Loops} & \textbf{Points/Side} & \textbf{GN} & \textbf{CG} & \textbf{GN} & \textbf{CG} \\
\midrule
4 & 16 & 6 & 14.7 & 6 & 12.3 \\
8 & 16 & 6 & 20.8 & 6 & 14.5 \\
16 & 16 & 6 & 28.7 & 6 & 15.3 \\
32 & 16 & 7 & 43.3  & 7 & 16.7 \\
64 & 16 & 6 & 69.2  & 6 & 16.7 \\
128 & 16 & 9 & 135.9  & 9 & 16.8 \\
\bottomrule
\end{tabular} }
\caption{Comparison of numerical scalability for an increasing number of subdomains when rotation is and is not used within the null space.}\label{tab:weak_null}
\end{table}

\section{Conclusion and Limitations}

This paper presents first results indicating the potential effectiveness of additive overlapping Schwarz preconditioners for a simple formulation of a pose graph SLAM problem.
As a test case, we consider a robot traversing a unit square trajectory with an increasing number of loops and odometry points along each side of the square.
In this setting, a significant improvement of the conditioning of the problem was observed with the use of the preconditioners, and excellent weak scalability is demonstrated with the two-level preconditioner using the GDSW coarse space.
The interpretation of a simplified SLAM problem as a finite element discretization presented in this paper further provides intuition for why preconditioning techniques developed for discretized PDEs may also be applicable to SLAM problems. This observation motivates the exploration of domain decomposition methods as scalable solvers for large-scale SLAM optimization problems.

While the results presented in this paper are encouraging, several limitations of the current study should be acknowledged.
The numerical experiments are conducted exclusively on a synthetic, highly regular test case, in which the robot repeatedly traverses a unit square trajectory.
This setup implies a perfectly uniform graph structure, regularly spaced loop closures, and Gaussian measurement noise—conditions that are rarely encountered in real-world deployments.
However, this test case was deliberately chosen to generate a controlled family of problem instances that allow systematic variation of the subdomain size and the number of subdomains.
This enables us to investigate the sensitivity of the method with respect to the subdomain size as well as its numerical scalability with respect to the overall problem size, in particular whether the iteration counts remain bounded independently of the problem size.

Furthermore, the current experiments place loop closure constraints within the overlap region, which may contribute to the favorable results.
Investigating other configurations will help clarify whether this choice is essential.

Future work will consider alternative nonlinear solvers such as Levenberg–Marquardt, as well as globalization and acceleration techniques in nonlinear optimization, including trust-region strategies and nonlinear domain decomposition methods.
We will further validate the approach on established real-world benchmarks such as the KITTI or EuRoC datasets
and study alternative decomposition strategies as well as different overlap sizes and their effect on the convergence behavior.

\section{Appendix}
\subsection{A Very Simple SLAM Problem} 
The purpose of this section is to show the relation between a simplified SLAM problem and classical finite element discretizations.  The model discussed here is deliberately elementary and serves as an explanatory tool; it is not identical to the nonlinear SLAM problem solved numerically in our experiments.

We consider a one-dimensional robot with scalar states $x_0, x_1, \dots, x_{10} \in \R$.

The robot moves from $x_0 = 0$ to $x_{10} = 10$ in ten steps of unit length (ground truth).  We assume that the odometry measurements are systematically biased.  Instead of the true step length $1$, the robot measures
\begin{equation*}
  \widetilde{\Phi}_{i,i+1} := 0.9, \qquad i = 0,\dots,9,
\end{equation*}
leading to the relations
\begin{equation*}
  x_{i+1} - x_i \approx 0.9.
\end{equation*}
The robot recognizes both its start and end positions exactly (end-closure and boundary conditions). Hence, we impose the hard constraints
\begin{equation*}
  x_0 = 0, \quad x_{10} = 10.
\end{equation*}
\subsection{Least-Squares Formulation}
The SLAM problem is formulated as a least-squares minimization problem.
For each odometry measurement, we define the residual
\begin{equation*}
  r_i(x) := (x_{i+1} - x_i) - \widetilde{\Phi}_{i,i+1} = (x_{i+1} - x_i) - 0.9.
\end{equation*}
The objective functional is
\begin{equation}
  \label{eq:1d_slam_functional}
  \mathcal{J}(x) = \frac12 \sum_{i=0}^{9} r_i(x)^{2},
\end{equation}
subject to the Dirichlet constraints on $x_0$ and $x_{10}$. Hence, our weight matrix $W$ is chosen as the identity matrix.

\subsection{Mechanical Interpretation and Equivalent Finite Element Formulation}
\label{sec:mech-fem}

Problem \eqref{eq:1d_slam_functional} is equivalent to a mechanical system consisting of a chain of linear axial bars: ten one-dimensional linear elastic bars; each bar has a rest length of $0.9$; node $0$ is fixed ($x_0 = 0$); node $10$ is prescribed to $x_{10} = 10$.

The total elastic energy of the system is defined as
\begin{equation*}
  \mathcal{E}(x) = \sum_{i=0}^{9} \frac12 \bigl((x_{i+1} - x_i) - 0.9\bigr)^2,
\end{equation*}
which coincides with the SLAM objective functional.
Thus, SLAM with end-closure can be interpreted as the computation of a mechanical equilibrium on a graph.

\subsection{Resulting Sparse Linear System}

The derivative of a single element energy
$
\frac12 \bigl((x_{i+1} - x_i) - 0.9\bigr)^2
$
yields the local stiffness matrix
\begin{equation*}
  K^{(i)} =
  \begin{pmatrix}
    1 & -1 \\
    -1 & 1
  \end{pmatrix},
\end{equation*}
acting on the degrees of freedom $(x_i,x_{i+1})$
and the local right hand side
\begin{equation*}
  b^{(i)} =
  \begin{pmatrix}
    -0.9\\
     0.9
  \end{pmatrix}.
\end{equation*}

Finite element assembly results in the global linear system
\begin{equation*}
  Kx = b
\end{equation*}
where $K$ is the $(11 \times 11)$ matrix
\begin{equation*}
  \label{eq:1d_laplace_matrix}
  K =
  \begin{pmatrix}
    1 & -1 & 0 & \cdots & 0 & 0 \\
    -1 & 2 & -1 & \ddots & \vdots & \vdots \\
    0 & -1 & 2 & \ddots & 0 & 0 \\
    \vdots & \ddots & \ddots & 2 & -1 & 0 \\
    0 & \cdots & 0 & -1 & 2 & -1 \\
    0 & \cdots & 0 & 0 & -1 & 1
  \end{pmatrix}.
\end{equation*}
and the right-hand side
$$
b=
\begin{pmatrix}
    -0.9, & 0, & \ldots & ,0, & 0.9
\end{pmatrix}^T .
$$
This is the classical tridiagonal matrix with stencil $(-1,\;2,\;-1)$, which discretizes the one-dimensional Laplace operator.


Finally, Dirichlet boundary conditions are imposed in the classical strong form by matrix modification: rows corresponding to fixed degrees of freedom are set to zero and the diagonal entries are then set to one. The prescribed values appear on the right-hand side.
For
\begin{equation*}
  x_0 = 0,
  \qquad
  x_{10} = 10,
\end{equation*}
this yields the modified system
\begin{equation*}
  \tilde K x = \tilde b
\end{equation*}
with
\begin{equation*}
  \tilde K =
  \begin{pmatrix}
    1 & 0 & 0 & \cdots & 0 & 0 \\[0.5ex]
    -1 & 2 & -1 & \ddots & 0 & 0 \\[0.5ex]
    0 & -1 & 2 & \ddots & 0 & 0 \\[0.5ex]
    \vdots & \ddots & \ddots & 2 & -1 & \vdots \\[0.5ex]
    0 & 0 & 0 & -1 & 2 & -1 \\[0.5ex]
    0 & 0 & 0 & \cdots & 0 & 1
  \end{pmatrix},
  \qquad
  \tilde b =
  \begin{pmatrix}
    0 \\[0.5ex] 0 \\[0.5ex] 0 \\[0.5ex] \vdots \\[0.5ex] 0 \\[0.5ex] 0 \\[0.5ex] 10
  \end{pmatrix}.
\end{equation*}

The resulting linear system is non-symmetric but symmetry can be reestablished, e.g., to allow for conjugate gradients, by an elimination step.  This removes the off-diagonal non-zero entries from the columns corresponding to the Dirichlet boundary.
We obtain
\begin{equation*}
  \hat K x = \hat b,
  \qquad
  x = (x_0,x_1,\dots,x_{10})^T,
\end{equation*}
with
\begin{equation*}
  \hat K =
  \begin{pmatrix}
    1 & 0 & 0 & \cdots & 0 & 0 & 0 \\[0.5ex]
    0 & 2 & -1 & \ddots & 0 & 0 & 0 \\[0.5ex]
    0 & -1 & 2 & \ddots & 0 & 0 & 0 \\[0.5ex]
    \vdots & \ddots & \ddots & \ddots & -1 & 0 & \vdots \\[0.5ex]
    0 & 0 & 0 & -1 & 2 & -1 & 0 \\[0.5ex]
    0 & 0 & 0 & 0 & -1 & 2 & 0 \\[0.5ex]
    0 & 0 & 0 & \cdots & 0 & 0 & 1
  \end{pmatrix},
  \qquad
  \hat b =
  \begin{pmatrix}
    0 \\[0.5ex]
    0 \\[0.5ex]
    0 \\[0.5ex]
    \vdots \\[0.5ex]
    0 \\[0.5ex]
    0 \\[0.5ex]
    10 \\[0.5ex]
    10
  \end{pmatrix}.
\end{equation*}

The solution of the system is
$x=(0,1,2,\ldots,10)^T.$
Mechanically, the excess length 
is distributed uniformly over all ten bars.

\subsection{Numerical Results} 
As a sanity check for our
overlapping Schwarz preconditioners, 
we present numerical tests for our very simple SLAM problem; see~Table~\ref{tab:laplace_scalability}.
As expected from the theory of domain decomposition methods,  only the two-level method is numerically scalable, i.e., the number of CG iterations is independent of the problem size.

\begin{table}[H]
\centering
\resizebox{.8\textwidth}{!}{
\begin{tabular}{cc||ccc||ccc}
\toprule
\multicolumn{2}{c||}{\textbf{}} & \multicolumn{3}{c||}{\textbf{One-Level Schwarz}} & \multicolumn{3}{c}{\textbf{GDSW}} \\
\cmidrule(lr){1-5}\cmidrule(lr){5-8}
\textbf{Subdomains} & \textbf{Bars} & \textbf{CG} & \textbf{$\lambda_{\rm min}$} & \textbf{$\lambda_{\rm max}$} & \textbf{CG} & \textbf{$\lambda_{\rm min}$} & \textbf{$\lambda_{\rm max}$} \\
\midrule
4 & 254 & 7 & $1.82\times 10^{-2}$ & 1.98 & 6 & 1.00 & 2.10 \\
8 & 506 & 16 & $4.79\times 10^{-3}$ & 1.99 & 6 & 1.00 & 2.10 \\
16 & 1010 & 32 & $1.22\times 10^{-3}$ & 2.00 & 6 & 1.00 & 2.10 \\
32 & 2018 & 65 & $3.05\times 10^{-4}$ & 2.00 & 6 & 1.00 & 2.10 \\
64 & 4034 & 132 & $7.64\times 10^{-5}$ & 2.00 & 6 & 1.00 & 2.10 \\
128 & 8066 & 262 & $1.91\times 10^{-5}$ & 2.00 & 6 & 1.00 & 2.10 \\
\bottomrule
\end{tabular}
}
\caption{Numerical scalability of one-level and two-level additive overlapping Schwarz preconditioners for varying subdomains and mesh resolution applied to the bar chain problem with overlap of 1. Only the two-level method is numerically scalable.}
\label{tab:laplace_scalability}
\end{table}

\subsection{Relation to the Full SLAM Problem}

The one-dimensional model discussed above is a linear scalar 
problem.
In contrast, the SLAM problem solved in the numerical experiments:
\begin{itemize}
\item involves vector-valued unknowns,
\item is nonlinear due to rotational components,
\item leads to state-dependent Jacobian matrices,
\item contains loop closure constraints inducing long-range couplings.
\end{itemize}

Nevertheless, the Gauss--Newton system matrices, see~\eqref{eq:slam-hessian}, arising in the full SLAM problem share essential properties with matrix discussed in the Appendix:
symmetry, positive definiteness after fixing the gauge, sparsity induced by local couplings, and a block structure corresponding to multi-component unknowns.

This analogy provides intuition for the use of domain decomposition methods in the SLAM-context and helps to explain the effectiveness of overlapping Schwarz preconditioners in this setting. It also motivates considering other types of DDM from structural mechanics~\cite{itersub_2009} as well.

\bibliographystyle{unsrt}
\bibliography{./all_refs}

\end{document}